\documentclass{amsart}
\usepackage[cp1251]{inputenc}
\usepackage[russian,english]{babel}
\usepackage{amscd}
\usepackage{amsmath}
\usepackage{graphicx}
\usepackage{amsfonts}
\usepackage{amssymb}
\allowdisplaybreaks

\newtheorem{theorem}{Theorem}[section]
\theoremstyle{plain}

\newtheorem{corollary}[theorem]{Corollary}

\newtheorem{definition}{Definition}[section]

\newtheorem{lemma}[theorem]{Lemma}

\newtheorem{proposition}[theorem]{Proposition}

\numberwithin{equation}{section}

\begin{document}
\author{Guy Roos}
\address{G.R.: 58 avenue de Versailles 75016 Paris France}
\email{guy.roos@laposte.net, roos@mx.amss.ac.cn}
\author{Weiping Yin}
\address{W.Y.: Dept. of Math., Capital Normal Univ., Beijing 100037, China}
\email{wyin@mail.cnu.edu.cn}
\thanks{The second author was supported in part by NSF of China (grant No: 19631010)
and by NSF of Beijing}
\title{New classes of domains with explicit Bergman kernel}
\date{}
\begin{abstract}We introduce two classes of egg type domains, built on general bounded
symmetric domains, for which we obtain the Bergman kernel in explicit formulas.
\end{abstract}
\maketitle

In 1921, S.~Bergman introduced a kernel function, which is now known as the
Bergman kernel function. It is well known that there exists a unique Bergman
kernel function for each bounded domain in $\mathbb{C}^{n}$. For which domains
can the Bergman kernel function be computed by explicit formulas? This is an
important problem. Explicit formulas of the Bergman kernel function can help
to solve important conjectures. We illustrate this point by two cases. Mostow
and Siu have given a counterexample to the important conjecture that the
universal covering of a compact K\"{a}hler manifold of negative sectional
curvature should be biholomorphic to the ball. In their counterexample the
explicit calculation of the Bergman kernel function and metric of the egg
domain
\[
\{z\in\mathbb{C}^{2}:|z_{1}|^{2}+|z_{2}|^{14}<1\}
\]
plays an essential role \cite{Y1}. Another example is Lu Qikeng conjecture: in
order to give a counterexample to the Lu Qikeng conjecture, an explicit
formula for the Bergman kernel function is used in \cite{Y2}. Therefore,
computation of the Bergman kernel function by explicit formulas is an
important research direction in several complex variables. Up to now, there
are still many mathematicians working in that direction.

If the group of holomorphic automorphisms of a bounded homogeneous domain is
known, then one can get its Bergman kernel function in explicit form. L.K.~Hua
\cite{Hua} has obtained the Bergman kernel functions in explicit form for the
four types of Cartan domains by using that method (it is called Hua's method).
It is well known that there are two exceptional Cartan domains, of dimension
16 and 27. Weiping Yin \cite{Y4} has obtained the Bergman kernel functions in
explicit form for these two exceptional Cartan domains by using the Hua's
method. A general expression of the Bergman kernel for all symmetric bounded
homogeneous domains can also be given using the theory of Jordan triple
systems; the Bergman kernel is then, up to a constant, a (negative) power of
the ''generic minimal polynomial'' (cf. Loos \cite{Lo77}).

For a non-symmetric homogeneous domain, we know that it is holomorphically
equivalent to a Siegel domain (or N-Siegel domain in the sense of Yichao Xu);
S.G.~Gindikin \cite{Y5} has computed the Bergman kernel functions in explicit
form for homogeneous Siegel domains by using generalized power functions.
Yichao Xu \cite{Y6} has also obtained the Bergman kernel functions in explicit
form for N-Siegel domains. In the middle of sixties, Jiaqing Zhong and Weiping
Yin \cite{Y7},\cite{Y8} have constructed some new types of non-symmetric
homogeneous domains and their extension spaces; Weiping Yin \cite{Y9}%
,\cite{Y10},\cite{Y11} has computed their Bergman kernel functions in explicit
form by using Hua's method. But their papers were published only at the
beginning of the eighties, after the end of Cultural Revolution in China.

Besides homogeneous domains, domains for which the Bergman kernel function can
be computed in explicit form are the egg domains (also called complex
ellipsoids or complex ovals). In general, an egg domain has the following
form:
\[
E(p_{1},\dots,p_{n})=\{(z_{1},\dots,z_{n})\in\mathbb{C}^{n}:|z_{1}|^{2p_{1}%
}+\dots+|z_{n}|^{2p_{n}}<1\},
\]
where $p_{1},\dots,p_{n}$ are positive real numbers. Here the $z_{j}$
$(j=1,\dots,n)$ are complex numbers. A more general case is obtained when the
$z_{j}$ are complex vectors, that is $z_{j}=\left(  z_{j1},\ldots,z_{jm_{j}%
}\right)  $ and $|z_{j}|^{2p_{j}}=[\sum_{k=1}^{m_{j}}|z_{jk}|^{2}]^{p_{j}}$;
the corresponding domain is then denoted by $E(p_{1},\dots,p_{n};m_{1}%
,\dots,m_{n})$ or simply by $E_{p;m}$. S.~Bergman (\cite{Y12}, p.82) has
computed the Bergman kernel function in explicit form on $E(p_{1},1)$ by
summing a series; although he states that $1/p_{1}$ is a positive integer, his
computation is valid for arbitrary $p_{1}$. The explicit form of the Bergman
kernel functions on $E(1,\ldots,1,p_{n})$ and on $E(1,\ldots,1,p_{n}%
;1,\ldots,1,m)$ were obtained by D'Angelo in \cite{Y13} and \cite{Y14}
respectively. In the case $1/p_{1},\dots,1/p_{n}$ are positive integers,
Zinov'ev \cite{Y15} obtained the explicit form of the Bergman kernel function
on $E(p_{1},\dots,p_{n})$. If $1/p_{1},\dots,1/p_{n-1}$ are positive integers
and $1/p_{n}=p$ is a positive real number, then the Bergman kernel function on
$E(p_{1},\dots,p_{n})$ can also be computed explicitly. If $p_{1},\dots,p_{n}$
are positive integers, then the Bergman kernel function on $E(p_{1}%
,\dots,p_{n})$ has an explicit expression in terms of multivariate
hypergeometric functions (\cite[Theorem 1]{Y16}).

We are not able to compute an explicit formula for the Bergman kernel function
on any egg domain. If there are two (or more) numbers among the positive real
numbers $1/p_{1},\dots,1/p_{n}$ which are not integers (nor inverse of
integers), then one cannot get an explicit formula for the Bergman kernel
function on $E(p_{1},\dots,p_{n})$ or on $E_{p;m}$. So, one needs to estimate
the Bergman kernel function on egg domains. This work has been done by Sheng
Gong and Xuean Zheng \cite{Y17},\cite{Y18}. In his thesis of doctoral degree
(Purdue Univ 1997-1998) titled ``The Bergman kernel on Reinhardt domains'',
Chieh-hsien Tiao has obtained an estimation for the Bergman kernel function on
general Reinhardt domains \cite{Ti99}.\bigskip\

During the period of the second author's stay at the Institut des Hautes
\'{E}tudes Scientifiques (IHES) in February 1998, he and G.~Roos introduced
the following four types of domains, called \emph{super-Cartan domains} or
\emph{Cartan-Hartogs domains} ($N,m,n$ are positive integers and $K>0$ is
real):
\begin{align*}
Y_{I}(N,m,n;K)  &  :=\{W\in\mathbb{C}^{N},Z\in R_{I}(m,n):|W|^{2K}%
<\det(I-ZZ^{T})\},\\
Y_{II}(N,p;K)  &  :=\{W\in\mathbb{C}^{N},Z\in R_{II}(p):|W|^{2K}<\det
(I-ZZ^{T})\},\\
Y_{III}(N,q;K)  &  :=\{W\in\mathbb{C}^{N},Z\in R_{III}(q):|W|^{2K}%
<\det(I-ZZ^{T})\},\\
Y_{IV}(N,n;K)  &  :=\{W\in\mathbb{C}^{N},Z\in R_{IV}(n):\\
&  \qquad|W|^{2K}<1-ZZ^{T}-[(ZZ^{T})^{2}-|ZZ^{\prime}|^{2}]^{1/2}\}.
\end{align*}
where $R_{I}(m,n),R_{II}(p),R_{III}(q)$ and $R_{IV}(n)$ denote respectively
the Cartan domains of first, second, third and fourth type in the sense of
L.K.~Hua. Here $Z^{T}$ denotes the conjugate and transpose of $Z$, $\det$
denotes the determinant of a square matrix.

Weiping Yin has got the Bergman kernel function in explicit form for these
four Cartan-Hartogs domains \cite{Yin99},\cite{Yin99b}. We would like to point
out that if we can get the Bergman kernel in explicit form for a domain, then
this domain is a good domain for research.

Quite recently, Weiping Yin has introduced the following four types of
domains, which may be called \emph{Cartan egg domains}:
\begin{align*}
CE_{I}(M,N,m,n;K)  &  :=\{W_{1}\in\mathbb{C}^{M},W_{2}\in\mathbb{C}^{N},Z\in
R_{I}(m,n):\\
&  \qquad|W_{1}|^{2}+|W_{2}|^{2K}<\det(I-ZZ^{T})\},\\
CE_{II}(M,N,p;K)  &  :=\{W_{1}\in\mathbb{C}^{M},W_{2}\in\mathbb{C}^{N},Z\in
R_{II}(p):\\
&  \qquad|W_{1}|^{2}+|W_{2}|^{2K}<\det(I-ZZ^{T})\},\\
CE_{III}(M,N,q;K)  &  :=\{W_{1}\in\mathbb{C}^{M},W_{2}\in\mathbb{C}^{N},Z\in
R_{III}(q):\\
&  \qquad|W_{1}|^{2}+|W_{2}|^{2K}<\det(I-ZZ^{T})\},\\
CE_{IV}(M,N,n;K)  &  :=\{W_{1}\in\mathbb{C}^{M},W_{2}\in\mathbb{C}^{N},Z\in
R_{IV}(n):\\
&  \qquad|W_{1}|^{2}+|W_{2}|^{2K}<1-ZZ^{T}-[(ZZ^{T})^{2}-|ZZ^{\prime}%
|^{2}]^{1/2}\}.
\end{align*}

More generally, we will consider in this paper the two following classes of
domains, built on an arbitrary irreducible bounded circled homogeneous domain
$\Omega$:
\begin{align*}
Y(q,\Omega;k)  &  :=\left\{  \left(  W,Z\right)  \in\mathbb{C}^{q}\times
\Omega;\left|  W\right|  ^{2k}<N(Z,Z)\right\}  ,\\
E(p,q,\Omega;k)  &  :=\left\{  \left(  W_{1},W_{2},Z\right)  \in\mathbb{C}%
^{p}\times\mathbb{C}^{q}\times\Omega;\ \left|  W_{1}\right|  ^{2}+\left|
W\right|  ^{2k}<N(Z,Z)\right\}  ,
\end{align*}
where $N(Z,Z)$ is the \emph{generic norm} of $\Omega$, which is the proper
generalization of \linebreak $\det(I-ZZ^{T})$.

If we are able to compute the Bergman kernel when $p=q=1$, then, by using
twice the \emph{Principle of inflation} (which will be stated below), one can
get the Bergman kernel functions of Cartan egg domains for general $p,q$.
Therefore, we consider first $Y(1,\Omega,k)$ and $E(1,1,\Omega;k)$. The domain
$E(1,1,\Omega;k)$ is not homogeneous, so we cannot use Hua's method to get the
Bergman kernel function in explicit form. Also the domain $E(1,1,\Omega;k)$ is
not Reinhardt, so we cannot get its Bergman kernel function by summing an
infinite series. Our new method is a combination of these two methods.\bigskip\

In Section \ref{SEC1}, we review the properties of Jordan triple systems,
which are necessary for the general definition of $Y(q,\Omega;k)$ and
$E(p,q,\Omega;k)$. Section \ref{SEC2} is devoted to auxiliary results: the
computation of $\int_{\Omega}N(x,x)^{s}\alpha^{n}$, which appears to be a
simple consequence of an integral computed by Selberg, the notion of
``semi-Reinhardt domains'' and description of complete orthonormal systems for
these domains, and the ``Principle of inflation'', which allows to derive the
Bergman kernel for $Y(q,\Omega;k)$ and $E(p,q,\Omega;k)$ from the special case
$p=q=1$. In Section \ref{SEC4}, we study automorphisms of the domains
$Y(q,\Omega;k)$ and $E(p,q,\Omega;k)$ and compute their Bergman kernels. The
last section gives tables stating specific results for all types (classical
and exceptional) of bounded symmetric domains.

\section{Bounded symmetric domains and Jordan triple systems\label{SEC1}}

Hereunder we give a review of properties of the Jordan triple structure
associated to a complex bounded symmetric domain (see \cite{Lo77}, \cite{Ro99}).

\subsection{Jordan triple system associated to a bounded symmetric domain}

Let $\Omega$ be an irreducible bounded circled homogeneous domain in a complex
vector space $V$. Let $K$ be the identity component of the (compact) Lie group
of (linear) automorphisms of $\Omega$ leaving $0$ fixed. Let $\omega$ be a
volume form on $V$, invariant by $K$ and by translations. Let $\mathcal{K}$ be
the Bergman kernel of $\Omega$ with respect to $\omega$, that is, the
reproducing kernel of the Hilbert space $H^{2}(\Omega,\omega)=\mathrm{Hol}%
(\Omega)\cap L^{2}(\Omega,\omega)$. The Bergman metric at $z\in\Omega$ is
defined by
\[
h_{z}(u,v)=\partial_{u}\overline{\partial}_{v}\log\mathcal{K}(z).
\]
The \emph{Jordan triple product} on $V$ is defined by
\[
h_{0}(\{uvw\},t)=\partial_{u}\overline{\partial}_{v}\partial_{w}%
\overline{\partial}_{t}\log\mathcal{K}(z)\left|  _{z=0}\right.  .
\]
The triple product $(x,y,z)\mapsto\{xyz\}$ is complex bilinear and symmetric
with respect to $(x,z)$, complex antilinear with respect to $y$. It satisfies
the\emph{\ Jordan identity}
\[
\{xy\{uvw\}\}-\{uv\{xyw\}\}=\{\{xyu\}vw\}-\{u\{vxy\}w\}.
\]
The space $V$ endowed with the triple product $\{xyz\}$ is called a
\emph{(Hermitian) Jordan triple system}. For $x,y,z\in V$, denote by $D(x,y)$
and $Q(x,z)$ the operators defined by
\[
\{xyz\}=D(x,y)z=Q(x,z)y.
\]
The Bergman metric at $0$ is related to $D$ by
\[
h_{0}(u,v)=\operatorname*{tr}D(u,v).
\]
A Jordan triple system is called \emph{Hermitian positive }if
$(u|v)=\operatorname*{tr}D(u,v)$ is positive definite. As the Bergman metric
of a bounded domain is always definite positive, the Jordan triple system
associated to a bounded symmetric domain is Hermitian positive.

The \emph{quadratic representation }
\[
Q:V\longrightarrow\operatorname*{End}{}_{\mathbb{R}}(V)
\]
is defined by $Q(x)y=\frac12\{xyx\}$. The following fundamental identity for
the quadratic representation is a consequence of the Jordan identity:
\[
Q(Q(x)y)=Q(x)Q(y)Q(x).
\]
The \emph{Bergman operator} $B$ is defined by
\[
B(x,y)=I-D(x,y)+Q(x)Q(y),
\]
where $I$ denotes the identity operator in $V$. It is also a consequence of
the Jordan identity that the following fundamental identity holds for the
Bergman operator:
\[
Q(B(x,y)z)=B(x,y)Q(z)B(y,x).
\]
The Bergman operator gets its name from the following property:
\[
h_{z}\left(  B(z,z)u,v\right)  =h_{0}(u,v)\quad(z\in\Omega;\ u,v\in V).
\]
If $\Phi\in\left(  \operatorname*{Aut}\Omega\right)  _{0}$, the identity
component of the automorphism group of $\Omega$, the relation
\begin{equation}
B(\Phi x,\Phi y)=\mathrm{d}\Phi(x)\circ B(x,y)\circ\mathrm{d}\Phi(y)^{*}
\label{K12}%
\end{equation}
holds for $x,y\in\Omega$, where $~^{*}$ denotes the adjoint with respect to
the hermitian metric $h_{0}$. As a consequence, the Bergman kernel of $\Omega$
is given by
\begin{equation}
\mathcal{K}(z)=\frac1{\operatorname*{vol}\Omega}\frac1{\det B(z,z)}.
\label{K1}%
\end{equation}

\subsection{Spectral theory}

An Hermitian positive Jordan triple system is always \emph{semi-simple}, that
is the direct sum of a finite family of simple subsystems with component-wise
triple product.

As the domain $\Omega$ is assumed to be irreducible, the associated Jordan
triple system $V$ is \emph{simple}, that is $V$ is not the direct sum of two
non trivial subsystems.

An \emph{automorphism} $f:V\rightarrow V$ of the Jordan triple system $V$ is a
complex linear isomorphism preserving the triple product~:
$f\{u,v,w\}=\{fu,fv,fw\}$. The automorphisms of $V$ form a group, denoted
Aut$\;V$, which is a compact Lie group; we will denote by $K$ its identity component.

An element $c\in V$ is called \emph{tripotent} if $\{ccc\}=2c$.If $c$ is a
tripotent, the operator $D(c,c)$ annihilates the polynomial $T(T-1)(T-2)$.

Let $c$ be a tripotent. The decomposition $V=V_{0}(c)\oplus V_{1}(c)\oplus
V_{2}(c) $, where $V_{j}(c)$ is the eigenspace $V_{j}(c)=\left\{  x\in
V\;;\;D(c,c)x=jx\right\}  $, is called the \emph{Peirce decomposition} of $V$
(with respect to the tripotent $c$).

Two tripotents $c_{1}$ and $c_{2}$ are called \emph{orthogonal} if
$D(c_{1},c_{2})=0$. If $c_{1}$ and $c_{2}$ are orthogonal tripotents, then
$D(c_{1},c_{1)}$ and $D(c_{2},c_{2})$ commute and $c_{1}+c_{2}$ is also a tripotent.

A non zero tripotent $c$ is called \emph{primitive} if it is not the sum of
non zero orthogonal tripotents. A tripotent $c$ is \emph{maximal} if there is
no non zero tripotent orthogonal to $c$. The set of maximal tripotents is
equal to the \emph{Shilov boundary} of the domain $\Omega$.

A \emph{frame} of $V$ is a maximal sequence $(c_{1},\ldots,c_{r})$ of pairwise
orthogonal primitive tripotents. The frames of $V$ form a manifold
$\mathcal{F}$, which is called the \emph{Satake-Furstenberg boundary} of
$\Omega$.

Let $\mathbf{c}=(c_{1},\ldots,c_{r})$ be a frame. For $0\leq i\leq j\leq r$,
let
\[
V_{ij}(\mathbf{c})=\left\{  x\in V\mid D(c_{k},c_{k})x=(\delta_{i}^{k}%
+\delta_{j}^{k})x,\;1\leq k\leq r\right\}  \text{~: }%
\]
the decomposition $V=\bigoplus_{0\leq i\leq j\leq r}V_{ij}(\mathbf{c})$ is
called the \emph{simultaneous Peirce decomposition} with respect to the frame
$\mathbf{c}$.

Let $V$ be a simple Hermitian positive Jordan triple system. Then there exist
frames for $V$. All frames have the same number of elements, which is the
\emph{rank} $r$ of $V$. The subspaces $V_{ij}=V_{ij}(\mathbf{c})$ of the
simultaneous Peirce decomposition have the following properties: $V_{00}=0$~;
$V_{ii}=\mathbb{C}e_{i}$ ($0<i$); all $V_{ij}$'s ($0<i<j$) have the same
dimension $a$; all $V_{0i}$'s ($0<i$) have the same dimension $b$.

The \emph{numerical invariants} of $V$ (or of $\Omega$) are the rank $r$ and
the two integers
\begin{gather*}
a=\dim V_{ij}\;\;(0<i<j)\text{,}\\
b=\dim V_{0i}\;\;(0<i)\text{.}%
\end{gather*}
The \emph{genus} of $V$ is the number $g$ defined by
\[
g=2+a(r-1)+b\text{.}%
\]
The HPJTS $V$ and the domain $\Omega$ are said to be of \emph{tube type} if
$b=0$.

Let $V$ be a simple Hermitian positive Jordan triple system. Then any $x\in V
$ can be written in a unique way
\[
x=\lambda_{1}c_{1}+\lambda_{2}c_{2}+\cdots+\lambda_{p}c_{p}\text{,}%
\]
where $\lambda_{1}>\lambda_{2}>\cdots>\lambda_{p}>0$ and $c_{1},c_{2}%
\ldots,c_{p}$ are pairwise orthogonal tripotents. The element $x$ is
\emph{regular} iff $p=r$; then $(c_{1},c_{2},\ldots,c_{r})$ is a frame of $V$.
The decomposition $x=\lambda_{1}c_{1}+\lambda_{2}c_{2}+\cdots+\lambda_{p}%
c_{p}$ is called the \emph{spectral decomposition} of $x$.

\subsection{The generic minimal polynomial}

Let $V$ be a Jordan triple system of rank $r$. There exist polynomials
$m_{1},\ldots,m_{r}$ on $V\times\overline{V}$, homogeneous of respective
bidegrees $(1,1),\ldots,(r,r)$, such that for each regular $x\in V$, the
polynomial
\[
m(T,x,y)=T^{r}-m_{1}(x,y)T^{r-1}+\cdots+(-1)^{r}m_{r}(x,y)
\]
satisfies
\[
m(T,x,x)=\prod_{i=1}^{r}(T-\lambda_{i}^{2})\text{, }%
\]
where $x=\lambda_{1}c_{1}+\lambda_{2}c_{2}+\cdots+\lambda_{r}c_{r}$ is the
spectral decomposition of $x$. Here $\overline{V}$ denotes the space $V$ with
the conjugate complex structure. The polynomial
\[
m(T,x,y)=T^{r}-m_{1}(x,y)T^{r-1}+\cdots+(-1)^{r}m_{r}(x,y)
\]
is called the \emph{generic minimal polynomial} of $V$ (at $(x,y)$). The
(inhomogeneous) polynomial $N:V\times\overline{V}\rightarrow\mathbb{C}$
defined by
\[
N(x,y)=m(1,x,y)
\]
is called the \emph{generic norm}. The following identities hold:
\begin{align}
\det B(x,y)  &  =N(x,y)^{g},\label{K2}\\
\operatorname*{tr}D(x,y)  &  =g\,m_{1}(x,y)\text{.}\nonumber
\end{align}

\subsection{The spectral norm}

Let $V$ be an HPJTS. The map $x\mapsto\lambda_{1}$, where $x=\lambda_{1}%
c_{1}+\lambda_{2}c_{2}+\cdots+\lambda_{p}c_{p}$ is the spectral decomposition
of $x$ ($\lambda_{1}>\lambda_{2}>\cdots>\lambda_{p}>0)$ is a norm on $V$,
called the \emph{spectral norm}. The \emph{bounded symmetric domain} $\Omega$
is the unit ball of $V$ for the \emph{spectral norm}. It is also characterized
by the set of polynomial inequalities
\[
\left.  \frac{\partial^{j}}{\partial T^{j}}m(T,x,x)\right|  _{T=1}%
>0,\qquad0\leq j\leq r-1.
\]

\begin{proposition}
\label{PR1}Let
\[
\Phi:\mathcal{F}\times\left\{  \lambda_{1}>\lambda_{2}>\cdots>\lambda
_{r}>0\right\}  \longrightarrow V_{\mathrm{reg}}%
\]
be defined by
\[
\Phi((c_{1},\ldots,c_{r}),(\lambda_{1},\ldots,\lambda_{r}))=\sum_{j=1}%
^{r}\lambda_{j}c_{j}\text{~; }%
\]
here $V_{\mathrm{reg}}$ is the open dense subset of regular elements of $V$.
Then $\Phi$ is a diffeomorphism; its restriction
\[
\Phi_{0}:\mathcal{F}\times\left\{  1>\lambda_{1}>\lambda_{2}>\cdots
>\lambda_{r}>0\right\}  \longrightarrow\Omega_{\mathrm{reg}}=\Omega\cap
V_{\mathrm{reg}}%
\]
is a diffeomorphism onto the open dense subset $\Omega_{\mathrm{reg}}$ of
regular elements of $\Omega$.
\end{proposition}

Let $\alpha$ be the K\"{a}hler form on $V$ associated to the Hermitian inner
product $m_{1}$:
\begin{equation}
\alpha=\frac{\mathrm{i}}{2\pi}\partial\overline{\partial}m_{1}\text{.}
\label{K8}%
\end{equation}
In the following, we endow $V$ with the volume form $\omega=\alpha^{n}$
($n=\dim_{\mathbb{C}}V$), so that the volume of the unit ball associated to
$m_{1}$ is equal to $1$. The following property is well-known (see e.g.,
\cite[(5.1.1)]{Ko99}):

\begin{proposition}
\label{PR2}Let $V$ be an irreducible HPJTS of dimension $n$, rank $r$ and
invariants $a$, $b$. The pull-back of the volume element $\alpha^{n}$ by
$\Phi$ is
\[
\Phi^{\ast}\alpha^{n}=\Theta\wedge\prod_{j=1}^{r}\lambda_{j}^{2b+1}%
\prod_{1\leq j<k\leq r}\left(  \lambda_{j}^{2}-\lambda_{k}^{2}\right)
^{a}\;\mathrm{d}\lambda_{1}\wedge\ldots\wedge\mathrm{d}\lambda_{r}\text{, }%
\]
where $a$ and $b$ are the numerical invariants of $V$ and $\Theta$ is a
$K$-invariant volume form on $\mathcal{F}$.
\end{proposition}

\section{Auxiliary results\label{SEC2}}

\subsection{The Selberg integral\label{SSEC21}}

Using an integral of Selberg, we compute in this subsection the integral
$\int_{\Omega}N(x,x)^{s}\alpha^{n}$, where $N(x,x)$ is the generic norm of the
irreducible bounded symmetric domain $\Omega$. The integral of Selberg was
used by Kor\'{a}nyi \cite{K82} for computing the volume of a general bounded
symmetric domain. The integral $\int_{\Omega}N(x,x)^{s}\alpha^{n}$ has been
calculated by Hua \cite{Hua} for the four series of classical domains.

\begin{proposition}
\label{TH1}Let $\operatorname*{vol}\Omega$ and $\operatorname*{vol}%
\mathcal{F}$ be the volumes of $\Omega$ and $\mathcal{F}$ w.r. to $\alpha^{n}$
and $\Theta$:
\[
\operatorname*{vol}\Omega=\int_{\Omega}\alpha^{n},\quad\operatorname*{vol}%
\mathcal{F}=\int_{\mathcal{F}}\Theta.
\]
Then, for $s\in\mathbb{C}$, $\operatorname*{Re}s>-1$
\begin{equation}
\int_{\Omega}N(x,x)^{s}\alpha^{n}=F(s)\operatorname*{vol}\mathcal{F}%
=\frac{F(s)}{F(0)}\operatorname*{vol}\Omega,\label{EQ2}%
\end{equation}
with
\begin{gather}
F(s)=\frac{1}{2^{r}r!}\prod_{j=1}^{r}\frac{\Gamma(b+1+(j-1)\frac{a}{2}%
)\Gamma(s+1+(j-1)\frac{a}{2})\Gamma(j\frac{a}{2}+1)}{\Gamma(s+b+2+(r+j-2)\frac
{a}{2})\Gamma(\frac{a}{2}+1)},\label{S6}\\
\frac{F(s)}{F(0)}=\prod_{j=1}^{r}\frac{\Gamma(s+1+(j-1)\frac{a}{2}%
)\Gamma(b+2+(r+j-2)\frac{a}{2})}{\Gamma(1+(j-1)\frac{a}{2})\Gamma
(s+b+2+(r+j-2)\frac{a}{2})}.\nonumber
\end{gather}
\end{proposition}%

\proof
According to Propositions \ref{PR1} and \ref{PR2}, we have
\[
F(s)=\idotsint\limits_{1>\lambda_{1}>\lambda_{2}>\cdots>\lambda_{r}>0}%
\prod_{j=1}^{r}(1-\lambda_{j}^{2})^{s}\prod_{j=1}^{r}\lambda_{j}^{2b+1}%
\prod_{1\leq j<k\leq r}\left(  \lambda_{j}^{2}-\lambda_{k}^{2}\right)
^{a}\;\mathrm{d}\lambda_{1}\wedge\ldots\wedge\mathrm{d}\lambda_{r}\text{ ; }%
\]
by the change of variables $t_{j}=\lambda_{j}^{2}$, we have
\[
F(s)=\frac1{2^{r}}\idotsint\limits_{1>\lambda_{1}>\lambda_{2}>\cdots
>\lambda_{r}>0}\prod_{j=1}^{r}(1-t_{j})^{s}\prod_{j=1}^{r}t_{j}^{b}%
\prod_{1\leq j<k\leq r}\left(  t_{j}-t_{k}\right)  ^{a}\;\mathrm{d}t_{1}%
\wedge\ldots\wedge\mathrm{d}t_{r}.
\]
Extending to integration over the cube $\left[  0,1\right]  ^{r}$, we also
have
\[
F(s)=\frac1{2^{r}r!}\int_{0}^{1}\cdots\int_{0}^{1}\prod_{j=1}^{r}(1-t_{j}%
)^{s}\prod_{j=1}^{r}t_{j}^{b}\prod_{1\leq j<k\leq r}\left|  t_{j}%
-t_{k}\right|  ^{a}\;\mathrm{d}t_{1}\wedge\ldots\wedge\mathrm{d}t_{r}.
\]

The above integral has been evaluated by Selberg \cite{S44}:

\begin{theorem}
[Selberg]\label{THSelberg}For $\operatorname*{Re}x>0$, $\operatorname*{Re}%
y>0$, $\operatorname*{Re}z>-\min\left(  \frac{1}{n},\frac{\operatorname*{Re}%
x}{n-1},\frac{\operatorname*{Re}y}{n-1}\right)  $, one has
\begin{multline*}
\int_{0}^{1}\cdots\int_{0}^{1}\prod_{j=1}^{n}t_{j}^{x-1}(1-t_{j})^{y-1}%
\prod_{1\leq j<k\leq n}\left|  t_{j}-t_{k}\right|  ^{2z}\;\mathrm{d}%
t_{1}\ldots\mathrm{d}t_{n}\\
=\prod_{j=1}^{n}\frac{\Gamma(x+(j-1)z)\Gamma(y+(j-1)z)\Gamma(jz+1)}%
{\Gamma(x+y+(n+j-2)z)\Gamma(z+1)}.
\end{multline*}
\end{theorem}

See \cite{A87} for a simple proof.

Applying this result for $n\leftarrow r$, $x\leftarrow b+1$, $y\leftarrow s+1
$, $z\leftarrow\frac a2$, we obtain the expression of $F(s)$:
\[
F(s)=\frac1{2^{r}r!}\prod_{j=1}^{r}\frac{\Gamma(b+1+(j-1)\frac a2)\Gamma
(s+1+(j-1)\frac a2)\Gamma(j\frac a2+1)}{\Gamma(s+b+2+(r+j-2)\frac
a2)\Gamma(\frac a2+1)},
\]
valid for $\operatorname*{Re}s>-1$. It follows that
\[
\frac{F(s)}{F(0)}=\prod_{j=1}^{r}\frac{\Gamma(s+1+(j-1)\frac a2)\Gamma
(b+2+(r+j-2)\frac a2)}{\Gamma(1+(j-1)\frac a2)\Gamma(s+b+2+(r+j-2)\frac a2)}.
\]

This proves Proposition \ref{TH1}.

More precisely, $\int_{\Omega}N(x,x)^{s}\alpha^{n}$ is a very simple rational
function of $s$:

\begin{theorem}
\label{TH2}Let $V$ be an irreducible HPJTS of dimension $n$. There exists a
polynomial $\chi$ of degree $n$, the zeroes of which are all negative integers
or half integers, such that
\[
\chi(s)\int_{\Omega}N(x,x)^{s}\alpha^{n}%
\]
is independent of $s$. The polynomial $\chi$ is equal to
\begin{equation}
\chi(s)=\prod_{j=1}^{r}\left(  s+1+(j-1)\frac{a}{2}\right)  _{1+b+\left(
r-j\right)  a}.\label{S9}%
\end{equation}
\end{theorem}

Here we use the classical notation $\left(  s\right)  _{k}$ for the polynomial
of degree $k$ (``Pochhammer polynomial''):
\[
\left(  s\right)  _{k}=\prod_{j=0}^{k-1}\left(  s+j\right)  =s(s+1)\cdots
(s+k-1).
\]%

\proof
It follows from (\ref{EQ2}) and (\ref{S6}) that $\int_{\Omega}N(x,x)^{s}%
\alpha^{n}$ is, up to a constant, equal to
\begin{equation}
\frac{1}{\chi(s)}=\prod_{j=1}^{r}\frac{\Gamma(s+1+(j-1)\frac{a}{2})}%
{\Gamma(s+b+2+(r+j-2)\frac{a}{2})}.\label{S8}%
\end{equation}
This may also be written (changing $j$ into $r+1-j$ in the denominator)%
\begin{align*}
\chi(s)  & =\frac{\prod_{j=1}^{r}\Gamma(s+b+2+(2r-j-1)\frac{a}{2})}%
{\prod_{j=1}^{r}\Gamma(s+1+(j-1)\frac{a}{2})}\\
& =\prod_{j=1}^{r}\left(  s+1+(j-1)\frac{a}{2}\right)  _{1+b+\left(
r-j\right)  a}.
\end{align*}%

\endproof

As a consequence of Theorem \ref{TH2}, we have for $s\in\mathbb{C}$,
$\operatorname*{Re}s>$-1
\begin{equation}
\int_{\Omega}N(x,x)^{s}\alpha^{n}=\frac{\chi(0)}{\chi(s)}\operatorname*{vol}%
\Omega.\label{S11}%
\end{equation}
Note also that the degree of $\chi$:%
\[
\sum_{j=1}^{r}\left(  1+b+(r-j)a\right)  =r+rb+\frac{r(r-1)}{2}a
\]
is equal to the dimension of $V$.

\subsection{Semi-Reinhardt domains\label{SemiR}}

\begin{definition}
A bounded domain $D$ in $\mathbb{C}^{m+n}$ is called a \emph{semi-Reinhardt
domain} if $0\in D$ and if
\[
(e^{i\theta_{1}}w_{1},\ldots,e^{i\theta_{m}}w_{m},e^{i\theta}z_{1}%
,\cdots,e^{i\theta}z_{n})\in D
\]
for all $(w,z)\in D$ and $\theta_{1},\ldots,\theta_{m},\theta$ real.
\end{definition}

That means, $D$ is ``Reinhardt'' w.r. to $w$ and is circular w.r. to $z$.
Obviously, a semi-Reinhardt domain is a circular domain, but the converse is
not true.

It is well-known that: if $D_{1}$ is a Reinhardt domain containing the origin
in $\mathbb{C}^{m}$, then $\{w_{1}^{j_{1}}\ldots w_{m}^{j_{m}}\}$ is a
complete orthogonal system for $D_{1}$; if $D_{2}$ is a circular domain
containing the origin in $\mathbb{C}^{n}$, then a complete orthonormal system
for $D_{2}$ is given by
\[
\left\{  P_{ki};~k\in\mathbb{N},~1\leq i\leq m_{k}=\binom{n+k-1}k\right\}  ,
\]
where, for any fixed $k$, $\left\{  P_{k1},P_{k2},\ldots,P_{km_{k}}\right\}  $
is an orthonormal basis of the space of homogeneous polynomials of degree $k$
in $z_{1},\ldots,z_{n}$. We give hereunder a generalization of these facts for
semi-Reinhardt domains.

Let $D$ be a semi-Reinhardt domain in $\mathbb{C}^{m+n}$. For each multi-index
$j\in\mathbb{N}^{m}$, consider the weighted scalar product on polynomials on
$\mathbb{C}^{n}$ defined by
\[
\left(  f\mid g\right)  ^{(j)}=\int_{D}\left|  w^{j}\right|  ^{2}%
f(z)\overline{g(z)}.
\]
As $D$ is circular w.r. to $z$, polynomials of different degrees are
orthogonal for this scalar product. For each $k\in\mathbb{N}$, choose a basis
of the space $\mathcal{P}_{k}$ of homogeneous polynomials of degree $k$ in
$z_{1},\ldots,z_{n}$, which is orthonormal w.r. to $\left(  ~\mid~\right)
^{(j)}$:
\[
\left\{  P_{ki}^{(j)}\right\}  _{1\leq i\leq m_{k}},
\]
where $m_{k}=\dim\mathcal{P}_{k}=\binom{n+k-1}k$. Then $\left\{  P_{ki}%
^{(j)};~k\in\mathbb{N},~1\leq i\leq m_{k}\right\}  $ is a complete orthonormal
system for the space of holomorphic functions such that
\[
\int_{D}\left|  w^{j}\right|  ^{2}\left|  f(z)\right|  ^{2}<\infty.
\]

\begin{theorem}
\label{TH3}Let $D$ be a semi-Reinhardt domain in $\mathbb{C}^{m+n}$. Then a
complete orthonormal system for $D$ is given by
\[
\left\{  w_{1}^{j_{1}}\dots w_{m}^{j_{m}}P_{ki}^{(j)}(z);~j=(j_{1},\dots
,j_{m})\in\mathbb{N}^{m},~k\in\mathbb{N},~1\leq i\leq m_{k}\right\}  .
\]
\end{theorem}%

\proof
Let $\Phi_{jki}(w,z)=w^{j}P_{ki}^{(j)}(z)$. As $D$ is semi-Reinhardt, it is
clear that the system $\left\{  \Phi_{jki}\right\}  $ is orthonormal.

Let $f(z,w)$ be a holomorphic, square-integrable function on $D$. Let
$b_{jki}=\int_{D}f(w,z)\overline{\Phi_{jki}(w,z)}$. It suffices to show that
\begin{equation}
\int_{D}\left|  f(w,z)\right|  ^{2}=\sum_{j,k,!\leq i\leq m_{k}}\left|
b_{jki}\right|  ^{2}. \label{R1}%
\end{equation}
Let $D_{z}=\left\{  w;~(z,w)\in D\right\}  $. Then $D_{z}$ is ``Reinhardt'',
$f(w,z)=\sum b_{j}(z)w^{j}$,
\[
b_{j}(z)=\frac{\int_{w\in D_{z}}f(w,z)\overline{w}^{j}}{\int_{w\in D_{z}%
}\left|  w^{j}\right|  ^{2}},
\]
and $\int_{w\in D_{z}}\left|  f(w,z)\right|  ^{2}=\sum_{j}\left|
b_{j}(z)\right|  ^{2}\int_{w\in D_{z}}\left|  w^{j}\right|  ^{2}$. The
functions $b_{j}(z)$ are holomorphic, for they may also be obtained locally as
Cauchy integrals of $f(z,w)$ on the Shilov boundary of a polydisc in
$\mathbb{C}_{w}^{m}$. Let $\Omega$ be the projection of $D$ on the second
factor $\mathbb{C}^{n}$; then $\Omega$ is circular, $b_{j}$ is
square-integrable on $\Omega$, $b_{j}=\sum_{k,i}c_{jki}P_{ki}^{(j)},$with
\[
c_{jki}=\int_{\Omega}\left|  w^{j}\right|  ^{2}b_{j}(z)\overline{P_{ki}^{(j)}%
}(z)=\int_{\Omega}f(w,z)\overline{w}^{j}\overline{P_{ki}^{(j)}}b_{jki}=b_{jki}%
\]
and $\int_{\Omega}\left|  w^{j}\right|  ^{2}\left|  b_{j}(z)\right|  ^{2}%
=\sum_{~k\in\mathbb{N},~1\leq i\leq m_{k}}\left|  b_{jki}\right|  ^{2}$. The
equality (\ref{R1}) then easily follows.
\endproof

\subsection{The principle of inflation\label{SSEC23}}

Let $\Omega$ be a bounded complete Hartogs domain in $\mathbb{C}^{n+1}$:
\[
\Omega=\left\{  (z,\zeta);~z\in D,~\zeta\in\mathbb{C},~|\zeta|^{2}%
<\phi(z)\right\}  ,
\]
where $\phi$ is a bounded, positive, continuous function on some bounded
domain $D$ in $\mathbb{C}^{n}$. Due to the circular symmetry w.r. to the
one-dimensional variable, the Bergman kernel of $\Omega$ can be written
\[
\mathcal{K}_{\Omega}(z,\zeta)=L(z,\left|  \zeta\right|  ^{2}).
\]
The ``Principle of inflation'', given in \cite{BFS}, allows to compute the
Bergman kernel of the ``inflated'' domain $G$, obtained from $\Omega$ when
replacing $\zeta$ by an $m$-dimensional variable $Z$ and $\left|
\zeta\right|  ^{2}$ by $\left\|  Z\right\|  ^{2}=|Z_{1}|^{2}+\dots+|Z_{m}%
|^{2}$.

\begin{theorem}
\label{INFL}(\cite{BFS}) Let $G$ be defined by
\[
G=\left\{  (z,Z);~z\in D,~Z\in\mathbb{C}^{m},~\left\|  Z\right\|  ^{2}%
<\phi(z)\right\}  .
\]
The spaces $\mathbb{C}^{n+1}$ and $\mathbb{C}^{n+m}$ being respectively
endowed with the translation invariant volume forms $\beta(z)\wedge
\frac{\mathrm{i}}{2\pi}\partial\overline{\partial}|\zeta|^{2}$ and
$\beta(z)\wedge\left(  \frac{\mathrm{i}}{2\pi}\partial\overline{\partial
}\left\|  Z\right\|  ^{2}\right)  ^{m}$, the Bergman kernel function
$\mathcal{K}_{G}$ of $G$ is
\[
K_{G}(z,Z)=\frac{1}{m!}\left.  \frac{\partial^{m-1}}{\partial r^{m-1}%
}L(z,r)\right|  _{r=\left\|  Z\right\|  ^{2}}.
\]
\end{theorem}

Trivial, but basic example: $n=0$, $\phi=1$; $\Omega$ is the unit disc, $G$
the unit hermitian ball,
\[
\mathcal{K}_{\Omega}(\zeta)=\frac1{\left(  1-\left|  \zeta\right|
^{2}\right)  ^{2}},\qquad\mathcal{K}_{G}(Z)=\frac1{\left(  1-\left\|
Z\right\|  ^{2}\right)  ^{m+1}}.
\]%

\proof
See \cite[Subsection 2.2]{BFS}.

\section{Egg domains built on bounded symmetric domains\label{SEC4}}

\subsection{The domains $Y(q,\Omega;k)$ and $E(p,q,\Omega;k)$}

Let $\Omega$ be any bounded irreducible circled homogeneous domain. For
$q\in\mathbb{N}$ and $k\in\mathbb{R}$, $k>0$, define the domain $Y(q,\Omega
;k)$ by
\[
Y(q,\Omega;k):=\left\{  \left(  W,Z\right)  \in\mathbb{C}^{q}\times
\Omega;\left\|  W\right\|  ^{2k}<N(Z,Z)\right\}  ;
\]
for $p,q\in\mathbb{N}$ and $k\in\mathbb{R}$, $k>0$, define the domain
$E(p,q,\Omega;k)$ by
\[
E(p,q,\Omega;k):=\left\{  \left(  W_{1},W_{2},Z\right)  \in\mathbb{C}%
^{p}\times\mathbb{C}^{q}\times\Omega;\ \left\|  W_{1}\right\|  ^{2}+\left\|
W_{2}\right\|  ^{2k}<N(Z,Z)\right\}  .
\]
Note that, as it is the case for $\Omega$, the boundary of $Y(q,\Omega;k)$ and
$E(p,q,\Omega;k)$ is not smooth when the rank of the symmetric domain $\Omega$
is greater than $1$.

For $\Omega$ belonging to one of the four classical series, an explicit form
of the Bergman kernel for $Y(q,\Omega;k)$ has been obtained by Weiping Yin
(\cite{Yin99}, \cite{Yin99b}).

The spaces $\mathbb{C}^{p}$ and $\mathbb{C}^{q}$ are equipped with their
canonical Hermitian structure and with the volume forms
\begin{align*}
\omega_{1}(W_{1})  &  =\left(  \frac{\mathrm{i}}{2\pi}\partial\overline
{\partial}\left\|  W_{1}\right\|  ^{2}\right)  ^{p},\\
\omega_{2}(W_{2})  &  =\left(  \frac{\mathrm{i}}{2\pi}\partial\overline
{\partial}\left\|  W_{2}\right\|  ^{2}\right)  ^{q},
\end{align*}
which give volume $1$ to the Hermitian unit balls. The Bergman kernel of
$Y(q,\Omega;k)$ (resp. $E(p,q,\Omega;k)$) is considered with respect to the
volume form $\omega_{2}(W)\wedge\omega(Z)$ (resp. $\omega_{1}(W_{1}%
)\wedge\omega_{2}(W_{2})\wedge\omega(Z)$), where $\omega=\alpha^{n}$ is
defined by (\ref{K8}).

\subsection{Automorphisms}

We consider here the case where $p=q=1$ and we write for short
\begin{align*}
Y(\Omega;k)  &  =Y(1,\Omega;k),\\
E(\Omega;k)  &  =E(1,1,\Omega;k).
\end{align*}

Let $\Phi\in\left(  \operatorname*{Aut}\Omega\right)  _{0}$, the identity
component of the automorphism group of $\Omega$, and let $J\Phi(Z)=\det
D\Phi(Z)$ be the complex Jacobian of $\Phi$ at $Z\in\Omega$. Let $Z_{0}$ be
the inverse image of $0$ under $\Phi$: $\Phi(Z_{0})=0$. From the relations
$B(\Phi(Z),\Phi(T))=\operatorname*{d}\Phi(Z)\circ B(Z,T)\circ\operatorname*{d}%
\Phi(T)^{*}$, $B(Z,0)=I$ and $\det B(Z,T)=N(Z,T)^{g}$, we obtain
\begin{equation}
1=J\Phi(Z)N(Z,Z_{0})^{g}\overline{J\Phi(Z_{0})}, \label{K15}%
\end{equation}
which shows that $N(Z,Z_{0})$ never vanishes when $Z\in\Omega$. In particular,
we have
\begin{equation}
J\Phi(Z)=\frac1{\overline{J\Phi(Z_{0})}N(Z,Z_{0})^{g}} \label{K7}%
\end{equation}
and
\begin{equation}
\left|  J\Phi(Z_{0})\right|  ^{2}=\frac1{N(Z_{0},Z_{0})^{g}} \label{K16}%
\end{equation}
Consider the holomorphic function $A(\quad,Z_{0})$ defined by
\begin{equation}
A(Z,Z_{0})=\frac{N(Z_{0},Z_{0})^{\frac12}}{N(Z,Z_{0})}; \label{K11}%
\end{equation}
up to multiplication by a complex number of modulus $1$, $J\Phi(Z)$ is then
the $g$-th power of $A(Z,Z_{0})$. As a bounded circled symmetric domain
$\Omega$ is always convex and $Z\mapsto A(Z,Z_{0})$ is a non-vanishing
holomorphic function on $\Omega$, the holomorphic function $Z\mapsto\ln
A(Z,Z_{0})$ is well-defined on $\Omega$ (say, with the initial condition that
$\ln A(Z_{0},Z_{0})$ is real); for each $\lambda\in\mathbb{C}$, we then define
the holomorphic function $Z\mapsto A(Z,Z_{0})^{\lambda} $ by $A(Z,Z_{0}%
)^{\lambda} =\exp\lambda\ln A(Z,Z_{0})$.

\begin{lemma}
Let $\Omega$ be a bounded irreducible circled homogeneous domain of genus $g$.
Then for each $\Phi\in\left(  \operatorname*{Aut}\Omega\right)  _{0}$ and for
complex numbers $\alpha_{1},\alpha_{2}$ of modulus $1$, the map $\Psi
:E(\Omega;k)\rightarrow E(\Omega;k)$ defined by
\begin{equation}
\Psi(W_{1},W_{2},Z)=(\alpha_{1}A(Z,Z_{0})W_{1},\alpha_{2}A(Z,Z_{0})^{\frac
{1}{k}}W_{2},\Phi(Z)),\label{K5}%
\end{equation}
where $Z_{0}=\Phi^{-1}(0)$, is a holomorphic automorphism of $E(\Omega;k)$.
Moreover, we have
\begin{equation}
\left|  J\Psi(W_{1},W_{2},Z)\right|  =\left|  A(Z,Z_{0})\right|  ^{1+\frac
{1}{k}+g}.\label{K3}%
\end{equation}
In particular,
\begin{equation}
\left|  J\Psi(W_{1},W_{2},Z_{0})\right|  ^{2}=N(Z_{0},Z_{0})^{-1-\frac{1}%
{k}-g}.\label{K6}%
\end{equation}
\end{lemma}%

\proof
From (\ref{K7}), (\ref{K16}), (\ref{K11}) we deduce that
\begin{equation}
\left|  J\Phi(Z)\right|  =\left|  A(Z,Z_{0})\right|  ^{g}. \label{K4}%
\end{equation}
From this, it is easily deduced that $\Psi$ is a holomorphic automorphism of
$E(\Omega;k)$. The differential of $\Psi$ has the form
\[
\mathrm{d}\Psi(W_{1},W_{2},Z)=\left(
\begin{array}
[c]{ccc}%
\alpha_{1}A(Z,Z_{0}) & 0 & 0\\
0 & \alpha_{2}(A(Z,Z_{0}))^{\frac1k} & 0\\
\ast & * & \mathrm{d}\Phi(Z)
\end{array}
\right)  ,
\]
which implies
\[
\left|  J\Psi(W_{1},W_{2},Z)\right|  =\left|  A(Z,Z_{0})\right|
^{1+\frac1k+g},
\]
that is (\ref{K3}); (\ref{K6}) then follows from $A(Z_{0},Z_{0})=N(Z_{0}%
,Z_{0})^{-\frac12}$.%
\endproof

As a corollary (or directly along the same lines), we have

\begin{corollary}
Let $\Omega$ be a bounded irreducible circled homogeneous domain of genus $g$.
Then for each $\Phi\in\left(  \operatorname*{Aut}\Omega\right)  _{0}$ and for
a complex number $\alpha$ of modulus $1$, the map $\Psi_{0}:Y(\Omega
;k)\rightarrow Y(\Omega;k)$ defined by
\[
\Psi_{0}(W,Z)=(\alpha A(Z,Z_{0})^{\frac{1}{k}}W,\Phi(Z)),
\]
where $\Phi(Z_{0})=0$, is a holomorphic automorphism of $Y(\Omega;k)$. We
have
\[
\left|  J\Psi_{0}(W,Z)\right|  =\left|  A(Z,Z_{0})\right|  ^{\frac{1}{k}+g}%
\]
and
\begin{equation}
\left|  J\Psi_{0}(W,Z_{0})\right|  ^{2}=N(Z_{0},Z_{0})^{-\frac{1}{k}%
-g}.\label{K6Y}%
\end{equation}
\end{corollary}

\subsection{The Bergman kernel of $Y(q,\Omega;k)$}

In view of the Principle of inflation, we consider first $Y(\Omega
;k)=Y(1,\Omega;k)$. Let $\mathcal{K}_{Y}$ denote the Bergman kernel of
$Y(\Omega;k)$. Assume $\Phi\in\left(  \operatorname*{Aut}\Omega\right)  _{0}$
is an automorphism of $\Omega$ which maps $Z\in D$ to $0$. Let $W^{*}%
=A(Z,Z)^{\frac1k}W$; then by the transformation law of the Bergman kernel and
(\ref{K6Y}), we have
\[
\mathcal{K}_{Y}(W,Z)=\mathcal{K}_{Y}(W^{*},0)\ N(Z,Z)^{-\frac1k-g}.
\]
So the determination of $\mathcal{K}_{Y}$ reduces to computing $\mathcal{K}%
_{Y}(W^{*},0)$.

The domain $Y(\Omega;k)$ is clearly semi-Reinhardt in $\mathbb{C}^{1+n}$;
applying Theorem \ref{TH3}, it has a complete orthonormal system of the form
\[
\left\{  W^{j}P_{\ell i}^{(j)}(Z);~j\in\mathbb{N},~k\in\mathbb{N},~1\leq i\leq
m_{\ell}=\binom{n+\ell-1}\ell\right\}  ,
\]
where the $P_{\ell i}^{(j)}$'s are polynomials of degree $\ell$ on $V$; in
particular, $P_{01}^{(j)}$ is the constant polynomial $a_{j}$ defined by
\[
\left|  a_{j}\right|  ^{2}\int_{Y(\Omega;k)}\left|  W^{j}\right|  ^{2}%
\omega_{2}(W)\wedge\omega(Z)=1.
\]
The Bergman kernel of $Y(\Omega;k)$ at $(W,0)$ is then
\[
\mathcal{K}_{Y}(W,0)=\sum\left|  W^{j}P_{\ell i}^{(j)}(0)\right|  ^{2}%
=\sum\left|  a_{j}\right|  ^{2}\left|  W^{j}\right|  ^{2}.
\]

\begin{lemma}%
\begin{equation}
\left|  a_{j}\right|  ^{2}=\frac{1}{\chi(0)\operatorname*{vol}\Omega}\left(
j+1\right)  \chi\left(  \frac{j+1}{k}\right)  .\label{K9Y}%
\end{equation}
\end{lemma}%

\proof
We have
\begin{align*}
\left|  a_{j}\right|  ^{-2}  &  =\int_{Y(\Omega;k)}\left|  W^{j}\right|
^{2}\omega_{2}(W)\wedge\omega(Z)\\
&  =\int_{Z\in\Omega}\left(  \int_{\left|  W\right|  ^{2k}<N(Z,Z)}\left|
W\right|  ^{2j}\omega_{2}(W)\right)  \omega(Z);
\end{align*}
as
\[
\int_{\left|  W\right|  ^{2k}<R^{2}}\left|  W\right|  ^{2j}\omega_{2}%
(W)=2\int_{0}^{R^{1/k}}r^{2j+1}dr=\frac{R^{2(j+1)/k}}{j+1},
\]
then, by Theorem \ref{TH2},
\[
\left|  a_{j}\right|  ^{-2}=\frac1{j+1}\int_{\Omega}N(Z,Z)^{(j+1)/k}%
\omega(Z)=\frac1{j+1}\frac{\chi(0)}{\chi((j+1)/k)}\operatorname*{vol}\Omega.
\]
This proves (\ref{K9Y}).%
\endproof

So we have
\[
\mathcal{K}_{Y}(W,0)=\sum\left|  a_{j}\right|  ^{2}\left|  W^{j}\right|
^{2}=\frac1{\chi(0)\operatorname*{vol}\Omega}\sum_{j=0}^{\infty}\left(
j+1\right)  \chi\left(  \frac{j+1}k\right)  \left|  W\right|  ^{2j}.
\]
Let $F$ be defined by
\[
F(X)=\sum_{j=0}^{\infty}\left(  \frac{j+1}k\right)  \chi\left(  \frac
{j+1}k\right)  X^{j};
\]
then
\[
\mathcal{K}_{Y}(W,0)=\frac k{\chi(0)\operatorname*{vol}\Omega}F\left(  \left|
W\right|  ^{2}\right)
\]
and
\begin{align*}
\mathcal{K}_{Y}(W,Z)  &  =\mathcal{K}_{Y}(W^{*},0)\ N(Z,Z)^{-\frac1k-g}\\
&  =\frac k{\chi(0)\operatorname*{vol}\Omega}F\left(  \frac{\left|  W\right|
^{2}}{N(Z,Z)^{1/k}}\right)  N(Z,Z)^{-\frac1k-g}.
\end{align*}
So we have proved

\begin{theorem}
The Bergman kernel of $Y(\Omega;k)=Y(1,\Omega;k)$ is
\begin{equation}
\mathcal{K}_{Y}(W,Z)=\frac{k}{\chi(0)\operatorname*{vol}\Omega}F\left(
\frac{\left|  W\right|  ^{2}}{N(Z,Z)^{1/k}}\right)  N(Z,Z)^{-\frac{1}{k}%
-g},\label{BY1}%
\end{equation}
where
\begin{equation}
F(X)=\sum_{j=0}^{\infty}\left(  \frac{j+1}{k}\right)  \chi\left(  \frac
{j+1}{k}\right)  X^{j}.\label{BY}%
\end{equation}
\end{theorem}

As $\left(  \frac{j+1}k\right)  \chi\left(  \frac{j+1}k\right)  $ is a
polynomial function of $j$, $F$ is a rational function of $X$ (a finite linear
combination of derivatives of $(1-X)^{-1}$). So $\mathcal{K}_{Y}(W,Z)$ is the
product of the Bergman kernel of $\Omega$:
\[
\mathcal{K}_{\Omega}(Z)=\frac1{\operatorname*{vol}\Omega}\frac1{N(Z,Z)^{g}}%
\]
by
\[
\frac k{\chi(0)}F\left(  \frac{\left|  W\right|  ^{2}}{N(Z,Z)^{1/k}}\right)
N(Z,Z)^{-\frac1k}=G\left(  \left|  W\right|  ^{2},N(Z,Z)^{1/k}\right)  ,
\]
where $G(X_{1},X_{2})$ is a rational function of $X_{1},X_{2}$.

The Bergman kernel of $Y(q,\Omega;k)$ is then obtained from $\mathcal{K}_{Y}$
by using the Principle of inflation (Theorem \ref{INFL}):

\begin{corollary}
The Bergman kernel of $Y(q,\Omega;k)$ is
\begin{equation}
\mathcal{K}_{Y(q,\Omega;k)}(W,Z)=\frac{1}{q!}\frac{k}{\chi
(0)\operatorname*{vol}\Omega}F^{(q-1)}\left(  \frac{\left\|  W\right\|  ^{2}%
}{N(Z,Z)^{1/k}}\right)  N(Z,Z)^{-\frac{q}{k}-g},\label{BYq}%
\end{equation}
where $F$ is the rational function defined by (\ref{BY}).
\end{corollary}

\subsection{The Bergman kernel of $E(p,q,\Omega;k)$}

Let $\mathcal{K}_{E}$ denote the Bergman kernel of $E(\Omega;k)=E(1,1,\Omega
;k)$. Assume $\Phi\in\left(  \operatorname*{Aut}\Omega\right)  _{0}$ is an
automorphism of $\Omega$ which maps $Z\in D$ to $0$. Let $W_{1}^{*}%
=A(Z,Z)W_{1}$, $W_{2}^{*}=A(Z,Z)^{\frac1k}W_{2}$; then by the transformation
law of the Bergman kernel and (\ref{K6}), we have
\[
\mathcal{K}_{E}(W_{1},W_{2},Z)=\mathcal{K}_{E}(W_{1}^{*},W_{2}^{*}%
,0)\ N(Z,Z)^{-1-\frac1k-g}.
\]
So the determination of $\mathcal{K}_{E}$ reduces to computing $\mathcal{K}%
_{E}(W_{1}^{*},W_{2}^{*},0)$.

The domain $E(\Omega;k)$ is semi-Reinhardt in $\mathbb{C}^{2+n}$; applying
Theorem \ref{TH3}, it has a complete orthonormal system of the form
\[
\left\{  W_{1}^{j_{1}}W_{2}^{j_{2}}P_{\ell i}^{(j_{1},j_{2})}(Z);~j_{1}%
,j_{2},\ell\in\mathbb{N},~1\leq i\leq m_{\ell}=\binom{n+\ell-1}\ell\right\}
,
\]
where the $P_{\ell i}^{(j_{1},j_{2})}$'s are polynomials of degree $\ell$; in
particular, $P_{01}^{(j_{1},j_{2})}$ is the constant polynomial $a_{j_{1}%
j_{2}}$ defined by
\[
\left|  a_{j_{1}j_{2}}\right|  ^{2}\int_{E(\Omega;k)}\left|  W_{1}^{j_{1}%
}W_{2}^{j_{2}}\right|  ^{2}\omega_{1}(W_{1})\wedge\omega_{2}(W_{2}%
)\wedge\omega(Z)=1.
\]
The Bergman kernel of $E(\Omega;k)$ at $(W_{1},W_{2},0)$ is then
\begin{align}
\mathcal{K}_{E}(W_{1},W_{2},0)  &  =\sum\left|  W_{1}^{j_{1}}W_{2}^{j_{2}%
}P_{\ell i}^{(j_{1},j_{2})}(0)\right|  ^{2}\label{K10}\\
&  =\sum\left|  a_{j_{1}j_{2}}\right|  ^{2}\left|  W_{1}^{j_{1}}\right|
^{2}\left|  W_{2}^{j_{2}}\right|  ^{2}.\nonumber
\end{align}

\begin{lemma}
Let $h=j_{1}+1+\frac{j_{2}+1}{k}$; then
\begin{equation}
\left|  a_{j_{1}j_{2}}\right|  ^{2}=\frac{k}{\chi(0)\operatorname*{vol}\Omega
}\frac{\Gamma\left(  h+1\right)  \chi(h)}{\Gamma(j_{1}+1)\Gamma\left(
\frac{j_{2}+1}{k}\right)  }.\label{K9}%
\end{equation}
\end{lemma}%

\proof
We have
\begin{align*}
\left|  a_{j_{1}j_{2}}\right|  ^{-2}  &  =\int_{E(\Omega;k)}\left|
W_{1}^{j_{1}}W_{2}^{j_{2}}\right|  ^{2}\omega_{1}(W_{1})\wedge\omega_{2}%
(W_{2})\wedge\omega(Z)\\
&  =\int_{Z\in\Omega}\left(  \int_{\left|  W_{1}\right|  ^{2}+\left|
W_{2}\right|  ^{2k}<N(Z,Z)}\left|  W_{1}^{j_{1}}W_{2}^{j_{2}}\right|
^{2}\omega_{1}(W_{1})\wedge\omega_{2}(W_{2})\right)  \omega(Z)
\end{align*}
and, by standard computations,
\begin{align*}
&  \int_{\left|  W_{1}\right|  ^{2}+\left|  W_{2}\right|  ^{2k}<R^{2}}\left|
W_{1}^{j_{1}}W_{2}^{j_{2}}\right|  ^{2}\omega_{1}(W_{1})\wedge\omega_{2}%
(W_{2})\\
&  =4\int_{r_{1}^{2}+r_{2}^{2k}<R^{2},\ r_{j}\geq0}r_{1}^{2j_{1}+1}%
r_{2}^{2j_{2}+1}\mathrm{d}r_{1}\mathrm{d}r_{2}\\
&  =\frac4k\int_{s_{1}^{2}+s_{2}^{2}<R^{2},\ s_{j}\geq0}s_{1}^{2j_{1}+1}%
s_{2}^{\frac{2\left(  j_{2}+1\right)  }k-1}\mathrm{d}s_{1}\mathrm{d}s_{2}\\
&  =\frac4k\int_{0}^{R}\rho^{2j_{1}+\frac{2\left(  j_{2}+1\right)  }%
k+1}\mathrm{d}\rho\ \int_{0}^{\pi/2}\cos^{2j_{1}+1}\theta\ \sin^{\frac
{2\left(  j_{2}+1\right)  }k-1}\theta\ \mathrm{d}\theta\\
&  =\frac4k\frac{R^{2\left(  j_{1}+1+\frac{j_{2}+1}k\right)  }}{2\left(
j_{1}+1+\frac{j_{2}+1}k\right)  }\frac{\Gamma(j_{1}+1)\Gamma\left(
\frac{j_{2}+1}k\right)  }{2\Gamma\left(  j_{1}+1+\frac{j_{2}+1}k\right)  }\\
&  =\frac1k\frac{\Gamma(j_{1}+1)\Gamma\left(  \frac{j_{2}+1}k\right)  }%
{\Gamma\left(  j_{1}+1+\frac{j_{2}+1}k+1\right)  }R^{2\left(  j_{1}%
+1+\frac{j_{2}+1}k\right)  }.
\end{align*}
Let $h=j_{1}+1+\frac{j_{2}+1}k$; then, by Theorem \ref{TH2},
\begin{align*}
\left|  a_{j_{1}j_{2}}\right|  ^{-2}  &  =\frac1k\frac{\Gamma(j_{1}%
+1)\Gamma\left(  \frac{j_{2}+1}k\right)  }{\Gamma\left(  h+1\right)  }%
\int_{\Omega}N(Z,Z)^{h}\omega(Z)\\
&  =\frac1k\frac{\Gamma(j_{1}+1)\Gamma\left(  \frac{j_{2}+1}k\right)  }%
{\Gamma\left(  h+1\right)  }\frac{\chi(0)}{\chi(h)}\operatorname*{vol}\Omega,
\end{align*}
which proves (\ref{K9}).%
\endproof

From (\ref{K9}) and (\ref{K10}), we obtain
\[
\mathcal{K}_{E}(W_{1},W_{2},0)=\frac1{\operatorname*{vol}\Omega}\frac
k{\chi(0)}\sum_{j_{1},j_{2}=0}^{\infty}\frac{\Gamma\left(  h+1\right)
\chi(h)}{\Gamma(j_{1}+1)\Gamma\left(  \frac{j_{2}+1}k\right)  }\left|
W_{1}^{j_{1}}\right|  ^{2}\left|  W_{2}^{j_{2}}\right|  ^{2}.
\]
Let $b_{j}$ ($1\leq j\leq n+2$) be the constants defined by
\[
h(h-1)\chi(h)=\sum_{j=1}^{n+2}b_{j}\left(  h+1\right)  _{j}=\sum_{j=1}%
^{n+2}b_{j}\frac{\Gamma(h+j+1)}{\Gamma(h+1)};
\]
so
\begin{align*}
\left|  a_{j_{1}j_{2}}\right|  ^{2}  &  =\frac1{\operatorname*{vol}\Omega
}\frac k{\chi(0)}\frac{\Gamma\left(  h+1\right)  \chi(h)}{\Gamma
(j_{1}+1)\Gamma\left(  \frac{j_{2}+1}k\right)  }\\
&  =\frac k{\chi(0)\operatorname*{vol}\Omega}\sum_{j=1}^{n+2}b_{j}\frac
{\Gamma(h+j+1)\Gamma(h-1)}{\Gamma(h+1)\Gamma(j_{1}+1)\Gamma\left(  \frac
{j_{2}+1}k\right)  }%
\end{align*}
and
\begin{multline*}
\mathcal{K}_{E}(W_{1},W_{2},0)=\frac1{\operatorname*{vol}\Omega}\frac
k{\chi(0)}\\
\sum_{j_{1},j_{2}=0}^{\infty}\sum_{j=1}^{n+2}b_{j}\frac{\Gamma(h+j+1)\Gamma
(h-1)}{\Gamma(h+1)\Gamma(j_{1}+1)\Gamma\left(  \frac{j_{2}+1}k\right)
}\left|  W_{1}^{j_{1}}\right|  ^{2}\left|  W_{2}^{j_{2}}\right|  ^{2}.
\end{multline*}
Let $s=(j_{2}+1)/k$, $t_{1}=|W_{1}|^{2}$, $t_{2}=|W_{2}|^{2}$, then we have
\[
\mathcal{K}_{E}(W_{1},W_{2},0)=\Lambda(t_{1},t_{2})
\]
where
\begin{equation}
\Lambda(t_{1},t_{2})=\frac1{\operatorname*{vol}\Omega}\frac k{\chi(0)}%
\sum_{j=1}^{n+2}b_{j}\sum_{j_{2}=0}^{\infty}\left(  \sum_{j_{1}=0}^{\infty
}\frac{\Gamma(h+j+1)\Gamma(j_{1}+s)}{\Gamma(h+1)\Gamma(j_{1}+1)\Gamma\left(
s\right)  }t_{1}^{j_{1}}\right)  t_{2}^{j_{2}}. \label{K14}%
\end{equation}
But
\begin{align*}
\sum_{j_{1}=0}^{\infty}\frac{\Gamma(h+j+1)\Gamma(j_{1}+s)}{\Gamma
(h+1)\Gamma(j_{1}+1)\Gamma\left(  s\right)  }t_{1}^{j_{1}}  &  =\sum_{j_{1}%
=0}^{\infty}\frac{\Gamma(j_{1}+j+2+s)\Gamma(j_{1}+s)}{\Gamma(j_{1}%
+2+s)\Gamma(j_{1}+1)\Gamma\left(  s\right)  }t_{1}^{j_{1}}\\
&  =\frac{\Gamma(j+2+s)}{\Gamma(s+2)}{}_{2}F_{1}(s+2+j,s;s+2;t_{1})
\end{align*}
where ${}_{2}F_{1}(a,b;c;z)$ denotes the Gauss hypergeometric function
\[
{}_{2}F_{1}(a,b;c;z)=\sum_{m=0}^{\infty}\frac{(a)_{m}\left(  b\right)  _{m}%
}{(c)_{m}}\frac{z^{m}}{\Gamma(m+1)}.
\]
By Euler's transformation formula (see \cite[(1.2.10)]{GR90}, or
\cite{Euler}):
\[
{}_{2}F_{1}(a,b;c;z)=(1-z)^{c-a-b}{}_{2}F_{1}(c-a,c-b;c;z),
\]
we get
\[
{}_{2}F_{1}(s+2+j,s;s+2;t_{1})=(1-t_{1})^{-s-j}{}_{2}F_{1}(-j,2;s+2;t_{1}),
\]
where ${}_{2}F_{1}(-j,2;s+2;t_{1})$ is in fact a polynomial of degree $j$:
\[
{}_{2}F_{1}(-j,2;s+2;t_{1})=\sum_{m=0}^{j}\frac{(-j)_{m}\left(  2\right)
_{m}}{(s+2)_{m}}\frac{t_{1}^{m}}{\Gamma(m+1)}.
\]
(It can be expressed with help of a Jacobi polynomial, but we will not use
this fact here.) This proves
\begin{align}
&  \sum_{j_{1}=0}^{\infty}\frac{\Gamma(h+j+1)\Gamma(j_{1}+s)}{\Gamma
(h+1)\Gamma(j_{1}+1)\Gamma\left(  s\right)  }t_{1}^{j_{1}}\label{K17}\\
&  =\Gamma(j+2+s)(1-t_{1})^{-s-j}\sum_{m=0}^{j}\frac{(-j)_{m}\left(  2\right)
_{m}}{\Gamma(s+2+m)}\frac{t_{1}^{m}}{\Gamma(m+1)}.\nonumber
\end{align}
From (\ref{K14}) and (\ref{K17}), we deduce
\begin{multline*}
\Lambda(t_{1},t_{2})=\frac1{\operatorname*{vol}\Omega}\frac k{\chi(0)}%
(1-t_{1})^{-1/k}\sum_{j=1}^{n+2}b_{j}(1-t_{1})^{-j}\\
\sum_{m=0}^{j}(-j)_{m}\left(  2\right)  _{m}\frac{t_{1}^{m}}{\Gamma(m+1)}%
\sum_{j_{2}=0}^{\infty}\frac{\Gamma(j+2+\frac{j_{2}+1}k)}{\Gamma(\frac
{j_{2}+1}k+2+m)}\left(  \frac{t_{2}}{\left(  1-t_{1}\right)  ^{1/k}}\right)
^{j_{2}}.
\end{multline*}
Consider now, for $0\leq m\leq j$, the function $H_{jm}$ defined by
\[
H_{jm}(\lambda)=\sum_{p=0}^{\infty}\frac{\Gamma(j+2+\frac{p+1}k)}{\Gamma
(\frac{p+1}k+2+m)}\lambda^{p}=\sum_{p=0}^{\infty}\left(  \frac{p+1}%
k+2+m\right)  _{j-m}\lambda^{p};
\]
as $\left(  \frac{p+1}k+2+m\right)  _{j-m}$ is a polynomial of $p$, $H_{jm}$
is a rational function of $\lambda$ (more precisely, a finite linear
combination of derivatives of the function $(1-\lambda)^{-1}$). So we finally obtain:

\begin{theorem}
The Bergman kernel $\mathcal{K}_{E}$ of $E(\Omega;k)$ is
\begin{equation}
\mathcal{K}_{E}(W_{1},W_{2},Z)=\Lambda\left(  \frac{\left|  W_{1}\right|
^{2}}{N(Z,Z)},\frac{\left|  W_{2}\right|  ^{2}}{N(Z,Z)^{1/k}}\right)
\ N(Z,Z)^{-1-\frac{1}{k}-g},\label{BE1}%
\end{equation}
where $\Lambda$ is defined by
\begin{gather}
\Lambda(t_{1},t_{2})=\frac{1}{\operatorname*{vol}\Omega}\frac{k}{\chi
(0)}(1-t_{1})^{-1/k}H\left(  t_{1},\frac{t_{2}}{\left(  1-t_{1}\right)
^{1/k}}\right)  ,\label{BE}\\
H(t_{1},\lambda)=\sum_{j=1}^{n+2}b_{j}(1-t_{1})^{-j}\sum_{m=0}^{j}%
(-j)_{m}\left(  2\right)  _{m}\frac{t_{1}^{m}}{\Gamma(m+1)}H_{jm}\left(
\lambda\right)  ,\nonumber\\
H_{jm}(\lambda)=\sum_{p=0}^{\infty}\left(  \frac{p+1}{k}+2+m\right)
_{j-m}\lambda^{p},\nonumber\\
h(h-1)\chi(h)=\sum_{j=1}^{n+2}b_{j}\left(  h+1\right)  _{j}.\nonumber
\end{gather}
\end{theorem}

We see that $\mathcal{K}_{E}(W_{1},W_{2},Z)$ is the product of the Bergman
kernel of $\Omega$:
\[
\mathcal{K}_{\Omega}(Z)=\frac1{\operatorname*{vol}\Omega}\frac1{N(Z,Z)^{g}}%
\]
by the function
\[
\frac k{\chi(0)}H\left(  \frac{\left|  W_{1}\right|  ^{2}}{N(Z,Z)}%
,\frac{\left|  W_{2}\right|  ^{2}}{\left(  N(Z,Z)-\left|  W_{1}\right|
^{2}\right)  ^{1/k}}\right)  N(Z,Z)^{-1-\frac1k},
\]
where $H(X_{1},X_{2})$ is a rational function of $X_{1},X_{2}$.

The Bergman kernel $\mathcal{K}_{E(p,q,\Omega;k)}$ of $E(p,q,\Omega;k)$ is
then obtained from $\mathcal{K}_{E}$ using twice the Principle of inflation
(Theorem \ref{INFL}):

\begin{corollary}
The Bergman kernel of $E(p,q,\Omega;k)$ is
\begin{multline}
\mathcal{K}_{E(p,q,\Omega;k)}(W_{1},W_{2},Z)\label{BEpq}\\
=\frac{1}{p!q!}\Lambda^{(p-1,q-1)}\left(  \frac{\left\|  W_{1}\right\|  ^{2}%
}{N(Z,Z)},\frac{\left\|  W_{2}\right\|  ^{2}}{N(Z,Z)^{1/k}}\right)
\ N(Z,Z)^{-p-\frac{q}{k}-g},\nonumber
\end{multline}
where
\[
\Lambda^{(p-1,q-1)}\left(  t_{1},t_{2}\right)  =\frac{\partial^{p+q-2}\Lambda
}{\partial t_{1}^{p-1}\partial t_{2}^{q-1}}\left(  t_{1},t_{2}\right)
\]
is the partial derivative of the function $\Lambda$ defined by (\ref{BE}).
\end{corollary}

\section{Tables\label{SEC3}}

The following examples exhaust the list of simple Hermitian positive Jordan
triple systems (see \cite{Lo75}). The HPJTS occurring in the four infinite
series $I_{p,q}$, $II_{n}$, $III_{n}$, $IV_{n}$ are called \emph{classical};
the two HPJTS of type $V$ and $VI$ are called \emph{exceptional}. There is
some overlapping between the classical series, due to a finite number of
isomorphisms in low dimension. We give hereunder for each type:

\begin{itemize}
\item  the definition of the space $V$, its Jordan triple product, and the
corresponding bounded circled homogeneous domain;

\item  the generic norm;

\item  the numerical invariants $r,a,b,g=2+a(r-1)+b$;

\item  the polynomial
\[
\chi(s)=\prod_{j=1}^{r}\left(  s+1+(j-1)\frac{a}{2}\right)  _{1+b+\left(
r-j\right)  a}.
\]
from which one can deduce the behavior of the function
\[
s\mapsto\int_{\Omega}N(x,x)^{s}\alpha^{n}%
\]
and the Bergman kernels of $Y(q,\Omega;k)$ and $E(p,q,\Omega;k)$.
\end{itemize}

\subsection{Type I$_{m,n}$ $\left(  1\leq m\leq n\right)  $}

$V=\mathcal{M}_{m,n}(\mathbb{C})$ (space of $m\times n$ matrices with complex
entries), endowed with the triple product
\[
\{xyz\}=x^{t}\overline{y}z+z^{t}\overline{y}x.
\]
The domain $\Omega$ is the set of $m\times n$ matrices $x$ such that
$I_{m}-x^{t}\overline{x}$ is definite positive. The generic minimal polynomial
is
\[
m(T,x,y)=\operatorname*{Det}(TI_{m}-x^{t}\overline{y}),
\]
where $\operatorname*{Det}$ is the usual determinant of square matrices. The
numerical invariants are $r=m$, $a=2$, $b=n-m$, $g=m+n$. These HPJTS are of
tube type only for $m=n$.

Applying (\ref{S9}), we obtain the polynomial $\chi$:%
\[
\chi(s)=\prod_{j=1}^{m}\left(  s+j\right)  _{m+n+1-2j},
\]
which may also be written%
\[
\chi(s)=\prod_{j=1}^{m}(s+j)_{n}.
\]

\subsection{Type II$_{n}$ $\left(  n\geq2\right)  $}

$V=\mathcal{A}_{n}(\mathbb{C})$ (space of $n\times n$ alternating matrices)
with the same triple product as for Type I. The domain $\Omega$ is the set of
$n\times n$ alternating matrices $x$ such that $I_{n}+x\overline{x}$ is
definite positive.

\subsubsection{Type II$_{2p}$ ($n=2p$ even)}

The generic minimal polynomial is here given by
\[
m(T,x,y)^{2}=\operatorname*{Det}(TI_{n}+x\overline{y}).
\]
The numerical invariants are $r=\frac{n}{2}=p$, $a=4$, $b=0$, $g=2\left(
n-1\right)  $; these HPJTS are of tube type.

The polynomial $\chi$ is
\[
\chi(s)=\prod_{j=1}^{p}\left(  s+2j-1\right)  _{1+4\left(  p-j\right)  }%
=\prod_{j=1}^{p}(s+2j-1)_{2p-1}.
\]

\subsubsection{Type II$_{2p+1}$ ($n=2p+1$ odd)}

The generic minimal polynomial is given by
\[
Tm(T,x,y)^{2}=\operatorname*{Det}(TI_{n}+x\overline{y}).
\]
The numerical invariants are $r=\left[  \frac n2\right]  =p$, $a=4$, $b=2$,
$g=2(n-1)$; these HPJTS are not of tube type.

The polynomial $\chi$ is%
\[
\chi(s)=\prod_{j=1}^{p}\left(  s+2j-1\right)  _{3+4\left(  p-j\right)  }%
=\prod_{j=1}^{p}(s+2j-1)_{2p+1}.
\]

\subsection{Type III$_{n}$ $\left(  n\geq1\right)  $}

$V=\mathcal{S}_{n}(\mathbb{C})$ (space of $n\times n$ symmetric matrices) with
the same triple product as for Type I. The domain $\Omega$ is the set of
$n\times n$ symmetric matrices $x$ such that $I_{n}-x\overline{x}$ is definite
positive. The generic minimal polynomial is
\[
m(T,x,y)=\operatorname*{Det}(TI_{n}-x\overline{y}).
\]
The numerical invariants are $r=n$, $a=1$, $b=0$, $g=n+1$. These HPJTS are of
tube type.

The polynomial $\chi$ is
\[
\chi(s)=\prod_{j=1}^{n}\left(  s+\frac{j+1}{2}\right)  _{1+n-j}.
\]

\subsection{Type IV$_{n}$ $\left(  n\neq2\right)  $}

$V=\mathbb{C}^{n}$ with the quadratic operator defined by
\[
Q(x)y=q(x,\overline{y})x-q(x)\overline{y},
\]
where $q(x)=\sum x_{i}^{2},$ $q(x,y)=2\sum x_{i}y_{i}$. The domain $\Omega$ is
the set of points $x\in\mathbb{C}^{n}$ such that
\[
1-q(x,\overline{x})+\left|  q(x)\right|  ^{2}>0,\qquad2-q(x,\overline{x})>0.
\]
The generic minimal polynomial is
\[
m(T,x,y)=T^{2}-q(x,\overline{y})+q(x)q(\overline{y}).
\]
The numerical invariants are $r=2$, $a=n-2$, $b=0$, $g=n$. These HPJTS are of
tube type.

The polynomial $\chi$ is
\[
\chi(s)=\left(  s+1\right)  _{n-1}\left(  s+\frac{n}{2}\right)  .
\]

\subsection{Type V}

$V=\mathcal{M}_{2,1}(\mathbb{O}_{\mathbb{C}})$, the subspace of $\mathcal{H}%
_{3}(\mathbb{O}_{\mathbb{C}})$ consisting in matrices of the form
\[
\left(
\begin{array}
[c]{ccc}%
0 & a_{3} & \tilde{a}_{2}\\
\tilde{a}_{3} & 0 & 0\\
a_{2} & 0 & 0
\end{array}
\right)
\]
with the same quadratic operator as for type VI (see below). Here $\tilde{a}$
denotes the Cayley conjugate of $a\in\mathbb{O}_{\mathbb{C}}$. The generic
minimal polynomial is
\[
m(T,x,y)=T^{2}-(x|y)T+(x^{\sharp}|y^{\sharp})\text{. }%
\]
The domain $\Omega$ is the ``exceptional domain of dimension $16$'' defined
by
\[
1-(x|x)+(x\mid x^{\sharp})>0,\qquad2-(x|x)>0.
\]
The numerical invariants are $r=2$, $a=6$, $b=4$, $g=12$. This HPJTS is not of
tube type.

The polynomial $\chi$ is
\[
\chi(s)=\prod_{j=1}^{r}\left(  s+1+(j-1)\frac{a}{2}\right)  _{1+b+\left(
r-j\right)  a}.
\]
\[
\chi(s)=\left(  s+1\right)  _{11}\left(  s+4\right)  _{5};
\]
it can also be written
\[
\chi(s)=(s+1)_{8}(s+4)_{8}.
\]

\subsection{\textbf{Type VI}}

$V=\mathcal{H}_{3}(\mathbb{O}_{\mathbb{C}})$, the space of $3\times3$ matrices
with entries in the space $\mathbb{O}_{\mathbb{C}}$ of octonions over
$\mathbb{C}$, which are Hermitian with respect to the Cayley conjugation; the
quadratic operator is defined by
\[
Q(x)y=(x|y)x-x^{\sharp}\times\overline{y},
\]
where $\times$ denotes the Freudenthal product, $x^{\sharp}$ the adjoint
matrix in $\mathcal{H}_{3}(\mathbb{O}_{\mathbb{C}})$ and $(x|y)$ the standard
Hermitian product in $\mathcal{H}_{3}(\mathbb{O}_{\mathbb{C}})$ (for details,
see \cite{Ro92}). The domain $\Omega$ is the ``exceptional domain of dimension
$27$'' defined by
\begin{gather*}
1-(x|x)+(x^{\sharp}|x^{\sharp})-\left|  \det x\right|  ^{2}>0,\\
3-2(x|x)+(x^{\sharp}|x^{\sharp})>0,\\
3-(x|x)>0.
\end{gather*}
The generic minimal polynomial is
\[
m(T,x,y)=T^{3}-(x|y)T^{2}+(x^{\sharp}|y^{\sharp})T-\det x\det\overline
{y}\text{, }%
\]
where $\det$ denotes the determinant in $\mathcal{H}_{3}(\mathbb{O}%
_{\mathbb{C}})$. The numerical invariants are $r=3$, $a=8$, $b=0$, $g=18$.
This HPJTS is of tube type.

The polynomial $\chi$ is
\begin{align*}
\chi(s)  & =(s+1)_{17}(s+5)_{9}(s+9)\\
& =(s+1)_{9}(s+5)_{9}(s+9)_{9}.
\end{align*}

\bigskip\

\end{document}